\documentclass{stylefile}          

\usepackage{amssymb}
\usepackage{bbold}
\usepackage{graphicx}
\usepackage{tikz}
\usetikzlibrary{arrows,intersections}
\usepackage{pgfplots}
\pgfplotsset{compat=1.3}
\usepackage{braket,amsfonts,amsopn}

\newcommand{\bR}{\mathbf{R}}
\newcommand{\bS}{\mathbf{S}}

\newcommand{\p}{\partial}

\DeclareMathOperator{\trace}{tr}

\theoremstyle{definition}
\newtheorem{experiment}{Experiment}[section]

\definecolor{dred}{RGB}{200,80,80}
\definecolor{dgreen}{RGB}{70,140,70}
\definecolor{dblue}{RGB}{80,80,200}

\title{Numerical Solution of the Simple Monge--Amp\`ere Equation with Non-convex Dirichlet Data on Non-convex Domains}
\headlinetitle{Monge--Amp\`ere equation with non-convex data and domains}

\lastnameone{Jensen}
\firstnameone{Max}
\nameshortone{M.~Jensen}
\addressone{Department of Mathematics, University of Sussex, Brighton BN1 9QF}
\countryone{England}
\emailone{m.jensen@sussex.ac.uk}

\abstract{The existence of a unique numerical solution of the semi-Lagrangian method for the simple Monge--Amp\`ere equation is known independently of the convexity of the domain or Dirichlet boundary data---when the Monge--Amp\`ere equation is posed as Bellman problem. However, the convergence to the viscosity solution has only been proved on strictly convex domains. In this paper we provide numerical evidence that convergence of numerical solutions is observed more generally without convexity assumptions. We illustrate how in the limit multi-valued functions may be approximated to satisfy the Dirichlet conditions on the boundary as well as local convexity in the interior of the domain.}



\researchsupported{}

\firstpage{1}

\begin{document}


\section{Introduction}

The paper is concerned with the numerical computation of solutions of the so-called simple Monge--Amp\`ere equation
\begin{subequations} \label{MA_intro}
\begin{alignat}{2}\label{MA}
\det(D^2 u) &= \Bigl( \frac{f}{2} \Bigr)^2 &&\qquad \mbox{in } \Omega, \\
u(x) &= g(x) &&\qquad \mbox{on } \partial\Omega, \label{BC}
\end{alignat}
\end{subequations}
where $\Omega \subset \bR^2$, $g \in C(\bR^2)$ and $f \in C(\Omega)$ is non-negative: $f \ge 0$. 

This raises immediately the question how the notion of solution for \eqref{MA_intro} should be defined. Beyond classical solutions, the generalisations in the Aleksandrov and viscosity sense provide settings of less smooth solutions. Any such definition in the literature imposes implicitly or explicitly convexity properties on $u$. 

If $\Omega$ is non-convex but it is known that there exists a locally convex subsolution which which attains the boundary data pointwise, then the works \cite{GS93,G98} give criteria for the existence of a unique solution. If $\Omega$ is convex (possibly not strictly) and $g$ non-convex then \cite{B} examines in the context of the closely related Gauss curvature problem the possibility of solutions which are in a certain sense multi-valued on $\partial \Omega$. We shall return to these results in sections \ref{sec:nonconvexdomain} and \ref{sec:nonconvexdata}.

To our knowledge there has been no systematic analysis of the well-posedness of \eqref{MA_intro} for the combination of a non-convex domain and general non-convex boundary data. It is therefore interesting if numerical methods can provide an insight into \eqref{MA_intro} in this case. The basis for our study is \cite{FJ} where a semi-Lagrangian numerical scheme is proposed, which uses Krylov's Bellman reformulation of \eqref{MA_intro}. While the convergence proof to viscosity solutions in \cite{FJ} requires strict convexity of $\Omega$ in order to make a comparison principle available, we highlight two results in that work which do not impose any form of convexity on $\Omega$ or $g$:

\begin{enumerate}[(a)]
\item For any finite element mesh there exists a unique numerical solution $u_i$, which is the limit of a globally converging semi-smooth Newton method.
\item There exists a constant $C$ (depending on $\Omega$, $f$ and $g$ only) such that $\| u_i \|_{L^\infty} \le C$ for all numerical solutions $u_i$ independently of the mesh.
\end{enumerate}

Thus at this point we know that as the finite element mesh size $h$ approaches $0$ some subsequences of numerical solutions converge weakly$^*$ to functions in $L^\infty$, which may be candidate solutions of \eqref{MA_intro} of some form.

The aims of this paper are
\begin{enumerate}[(a)]
\item to provide numerical evidence that indeed not just subsequences but the whole sequence of numerical solutions converges, provided that the stencil size is scaled appropriately with the mesh size;
\item to highlight that the Bakelman interpretation of solutions might extend to the non-convex setting, consistently with the findings in \cite{GS93,G98};
\item to study the performance of the numerical scheme in the non-convex setting, e.g.~the robustness of the semi-smooth Newton solver;
\item to present computations on general domains, e.g.~which are not Lipschitz.
\end{enumerate}

\section{The Bellman formulation and its semi-Lagrangian approximation}

A difficulty when solving \eqref{MA_intro} is that the Monge--Amp\`ere operator $v \mapsto \det(D^2 u)$ is only elliptic on the set of convex functions. It is therefore convenient to reformulate the problem in such a way that the set of solutions remains unchanged but ellipticity is established on the whole function space.

For that purpose we define the Bellman operator
\begin{alignat}{2}\label{HJB_operator}
H(A,f) &:=\sup_{B \in \bS_1} \Bigl( -B:A + f \sqrt[d]{\det B} \Bigr) 
&&\qquad \forall A\in \bS, f \in [0,\infty),
\end{alignat}
where $\bS$ is the set of symmetric $2 \times 2$ matrices, $\bS_+ := \{ A \in \bS;\, A \ge 0 \}$ and $\bS_1 := \{ B \in \bS_+;\, \trace B = 1 \}$, and consider the boundary value problem
\begin{subequations} \label{EHJB}
\begin{align} \label{EHJBeq}
H\bigl(D^2u(x),f(x)\bigr) & = 0 && \forall x\in \Omega,\\[1mm]
u(x) & = g(x) && \forall x\in \p\Omega. \label{EHJBbc}
\end{align}
\end{subequations}
It was shown \cite{FJ} that
\begin{align} \label{eq:equivalence}
& \{ v \in C(\overline{\Omega}) : \text{$v$ is convex and a viscosity solution of \eqref{MA}} \}\\
= \; & \{ v \in C(\overline{\Omega}) : \text{$v$ is a viscosity solution of \eqref{EHJBeq}} \}. \nonumber
\end{align}
This statement does not require boundedness or convexity of $\Omega$ and does in this form not refer to the boundary conditions.

To discretize \eqref{EHJBeq} we write $B = \sum_i \lambda_i \, e_i \, e_i^{\sf T} \in \bS_1$ with the eigenvalues $\lambda_i$ and normalised eigenvectors $e_i$, so that for smooth functions $\phi$
\begin{align*}
    B : D^2 \phi(x) & = \sum_i \lambda_i (e_i \, e_i^{\sf T} ) : D^2 \phi(x) = \sum_i \lambda_i \, \partial_{e_i,e_i}^2 \phi(x)\\
   &  = \sum_i \lambda_i \frac{\phi(x + k e_i) - 2 \phi(x) + \phi(x - k e_i)}{k^2} + \mathcal{O}(k^2).
\end{align*} 
Evaluating these central differences on a P1 finite element space at interior nodes, combined with nodal interpolation of $g$ on $\partial \Omega$, gives the numerical scheme of \cite{FJ}.

\begin{remark} \label{rem:classicalBC}
A comparison principle on the set of semi-continuous functions in respect to $H$ holds when classical boundary conditions are considered \cite{FJ}, where $\Omega$ may be non-convex. For the proof of convergence this comparison principle is applied to the upper and lower semi-continuous envelopes of the sequence of numerical solutions---in the setting of strictly convex domains where it is ensured that the envelopes attain the boundary conditions in the classical sense.

The numerical experiments in this paper highlight that on non-convex domains the envelopes may only satisfy the boundary conditions in a generalised form. It would be most convenient if the convergence proof could be translated to viscosity boundary conditions as they are defined either in \cite{BS} or alternatively in the form of \cite{UG}. However, in \cite{JS} counterexamples were given to show that comparison principles with these kinds of viscosity boundary conditions do in general not hold on the spaces of semi-continuous functions (even if $\Omega$ is convex), screening out trivial extensions of the proof.
\end{remark}

\section{The Monge--Amp\`ere equation on non-convex domains} \label{sec:nonconvexdomain}

In the light of Remark \ref{rem:classicalBC}, we wish to identify settings in which the exact solution of \eqref{MA_intro} satisfies the boundary conditions classically. If was shown in \cite{GS93} that this is guaranteed provided one knows of the existence of a strict subsolution and sufficient smoothness of data and domain. Indeed \cite{GS93} covers a more general equation with dependence on first-order derivatives.

\begin{theorem}[Guan and Spruck `93] \label{thm:GS93}
Let $\Omega \subset \bR^n$ be a smooth bounded domain with boundary components $\partial \Omega = (\gamma_1, \ldots, \gamma_m)$. Assume that there is a smooth, strictly locally convex function $\underline{u}$ in $\overline{\Omega}$ satisfying
\begin{align*}
\det(\underline{u}_{ij}) & \ge \psi(x, \underline{u}, \nabla \underline{u}) + \delta_0 \quad \text{in } \Omega,\\
\underline{u} & = \phi \quad \text{on } \partial \Omega,
\end{align*}
where $\phi$ and $\psi$ are smooth, $\psi > 0$ and $\psi^{1/n}(x,u,p)$ is convex in $p$. Then there is a smooth locally convex solution $u \in C^\infty(\overline{\Omega})$ to the Monge--Amp\`ere boundary-value problem:
\begin{align*}
\det(u_{ij}) & = \psi(x, u, \nabla u) \quad \text{in } \Omega,\\
u & = \phi \quad \text{on } \partial \Omega.
\end{align*}
If $\psi_u \ge 0$, then the solution is unique. 
\end{theorem}

A similar result can be found in \cite{G98}, with vanishing right-hand side but a more precise specification of the required regularity.

\begin{theorem}[Guan `98] \label{thm:G98}
Assume that $\partial \Omega$ is in $C^{3,1}$ and $\phi \in C^{3,1}(\partial \Omega)$. Suppose there exists a locally strictly convex function $\underline{u} \in C^2(\overline{\Omega})$ with $\underline{u} = \phi$ on $\partial \Omega$. There there is a unique locally convex weak solution of
\[
\det(u_{ij}) = 0 \quad \text{in} \quad \Omega, \quad u = \phi \quad \text{on} \quad \partial \Omega
\]
in $C^{1,1}(\overline{\Omega})$.
\end{theorem}

Our first computational experiment concerns the convergence of the numerical approximations to an exact solution, for which we know from the above that the boundary conditions are admitted classically.

\begin{experiment}[Quartic problem on L-shape] \label{exp:quarticLshape}
We approximate the exact solution $u(x) = (x_1^2 + x_2^2)^2$ on the L-shape:
\[
\Omega = \bigl[ (0,1) \times (-1,1) \bigr] \cup \bigl[ (-1,1) \times (0,1) \bigr].
\]
The quasi-uniform grid has at the coarsest level $58$ nodes and at the finest level after $6$ uniform refinements $168,961$ nodes. The stencil diameter $k$ is, away from the boundary, represented through $k = m \cdot h$ by a fixed positive factor $m$ and the (average) mesh size $h$. Near the boundary, so where $m \cdot h$ is larger than the distance to $\partial \Omega$, the stencil is reduced in size to remain within $\Omega$, cf.~\cite{FJ}.

\begin{figure}[t]
\begin{center}
\begin{tikzpicture}
    \begin{axis}[
        xlabel={Quartic problem on L-shape: Number of refinements},
        ylabel={$\| u - u_h \|_{\infty} / \| u \|_{\infty}$},
        ymode=log,
        legend entries={$m=2$,$m=4$,$m=8$,$m=16$,$m=32$,$m=64$},
        legend pos=south west,
        x=1.2cm,y=1.0cm,
        ymin=0.0007, ymax=0.1,
    ]

    \addplot plot coordinates { (0,  0.06281535) (1,  0.03476979) (2,  0.01827319) (3,  0.01344558) (4,  0.01171854) (5,  0.01099218) (6,  0.01066449) };
    \addplot plot coordinates { (0,  0.06867024) (1,  0.05886318) (2,  0.02540731) (3,  0.00965045) (4,  0.00491577) (5,  0.00359392) (6,  0.00320195) };
    \addplot plot coordinates { (0,  0.06963096) (1,  0.06465305) (2,  0.05483567) (3,  0.02251409) (4,  0.00714082) (5,  0.00253235) (6,  0.00130613) };
    \addplot plot coordinates { (0,  0.0670038 ) (1,  0.05887076) (2,  0.05929836) (3,  0.05252462) (4,  0.0215536 ) (5,  0.00643602) (6,  0.00188738) };
    \addplot plot coordinates { (0,  0.06667522) (1,  0.05768082) (2,  0.05840153) (3,  0.05724806) (4,  0.05184097) (5,  0.02127453) (6,  0.00623872) };
    \addplot plot coordinates { (0,  0.06742578) (1,  0.06010477) (2,  0.05175829) (3,  0.05780147) (4,  0.05670522) (5,  0.05165922) (6,  0.02119706) };
    \end{axis}
\end{tikzpicture}
\end{center}
\caption{Relative $L^\infty$-error for the test problem of Experiment \ref{exp:quarticLshape}.}
\label{fig:quarticLshape:rates}
\end{figure}
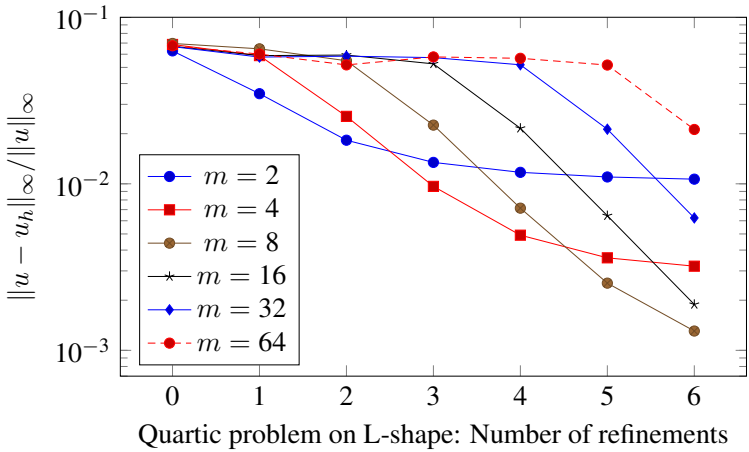

Figure \ref{fig:quarticLshape:rates} shows the decay of the relative $L^\infty(\Omega)$ approximation error for different choices of the multiplier $m$. Importantly, this gives numerical evidence of the convergence of the scheme on a non-convex domain. The table shows the multiplier $m$ which achieves the smallest relative $L^\infty$ error for a given number of degrees of freedom, and the number of Newton iterations of the respective computations to obtain a Newton step size of less than 5e-8 in the $\infty$-norm:

\begin{center}
\begin{tabular}{r||r|r|r|r|r|r|r}
DoFs & 58 & 197 & 721 & 2753 & 10753 & 42497 & 168961\\ \hline
$m$ & 2 & 2 & 2 & 4 & 4 & 8 & 8\\ 
rel.~$L^\infty$-error & 6.2e-2 & 3.5e-2 & 1.8e-2 & 9.6e-3 & 4.9e-3 & 2.5e-03 & 1.3e-3\\
\# Newton & 5 & 6 & 5 & 7 & 7 & 8 & 8
\end{tabular}
\end{center}

Overall the largest number of Newton iterations in this computational experiment is $9$, which occurs on the finest mesh with $m = 32$.
\end{experiment}

\section{The simple Monge--Amp\`ere equation with non-convex boundary data} \label{sec:nonconvexdata}

In computational experiments where the domain $\Omega$ is not strictly convex and there is no subsolution $\underline{u}$ as in the above Theorems \ref{thm:GS93} and \ref{thm:G98} we observe numerical solutions which appear to approximate multi-valued boundary data.

\begin{figure}[t]
\begin{center}
\includegraphics[width=9cm]{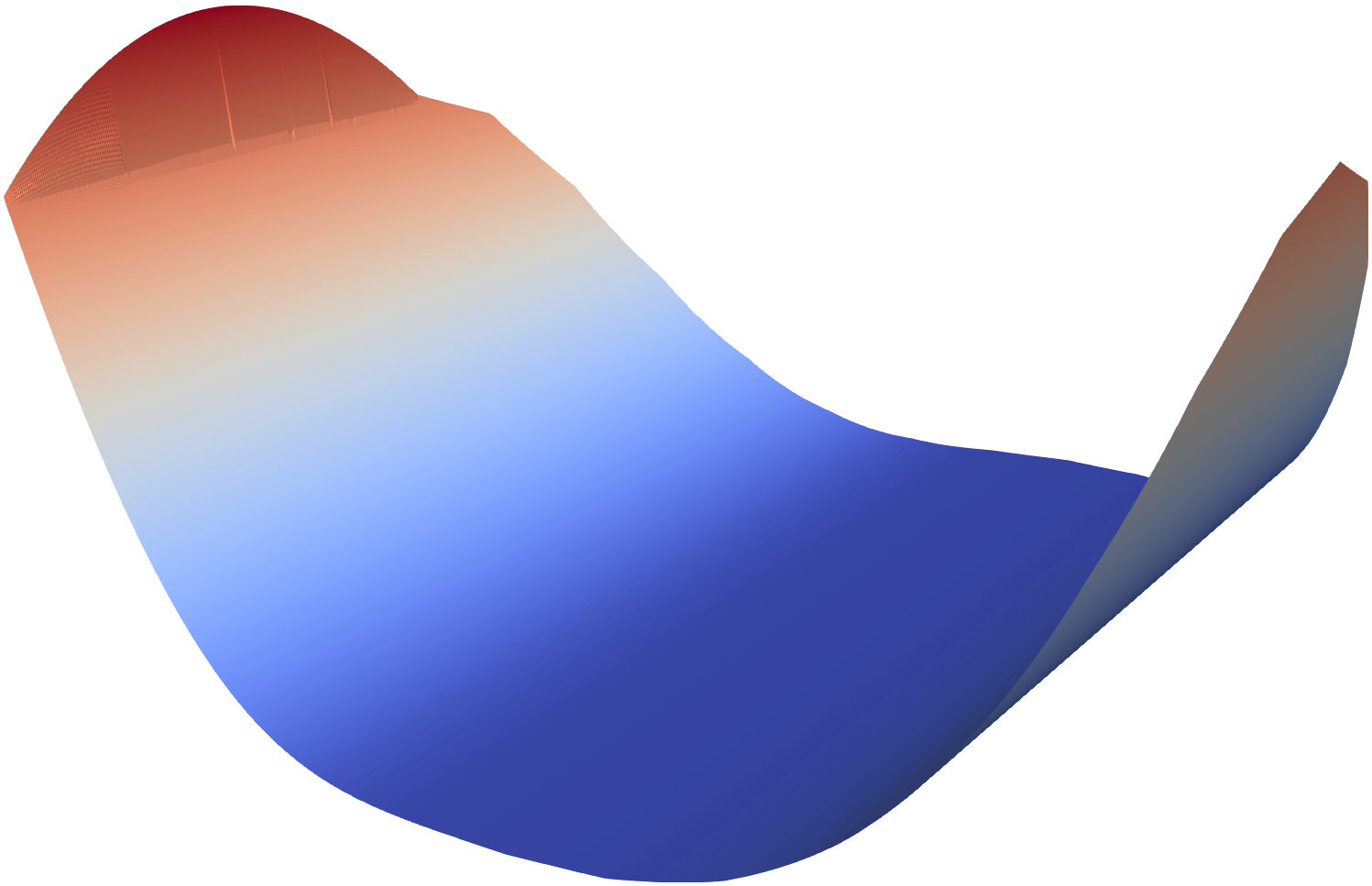}\\(a)

\bigskip

\includegraphics[width=9cm]{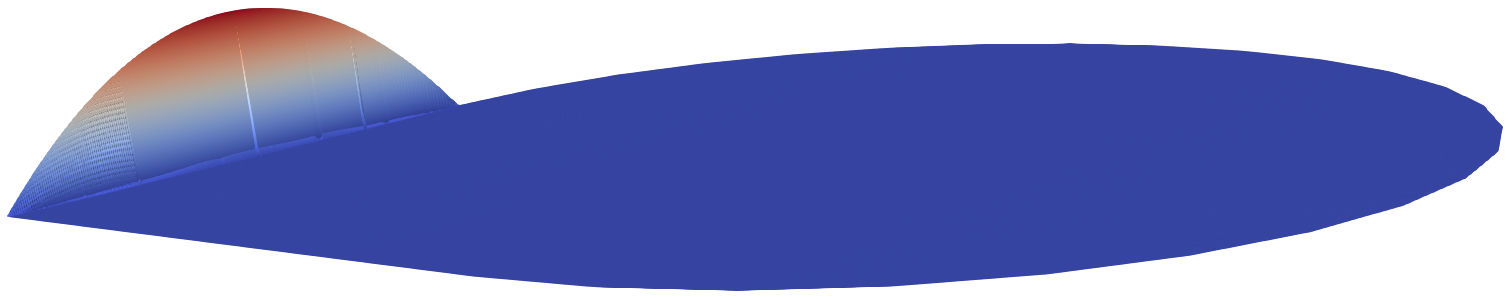}\\(b)
\end{center}
\caption{Plot (a) shows the numerical solution of Experiment \ref{exp:Bakelman} on the finest mesh. Notice the approximation to a multi-valued solution on the far side. Here the computational domain is given by the approximation of $\Omega$ with the coarsest mesh, which is not changed in the course of the refinement. In (b) we see the difference between the numerical solution and $x_1^4 - 1$, highlighting how localised the effect of the non-convex section of~$g$ on the numerical solution is.}
\label{fig:Bakelman:plot}
\end{figure}

\begin{figure}[t]
\begin{center}
\begin{tikzpicture}
    \begin{axis}[
        xlabel={Non-convex boundary data on convex domain: Number of refinements},
        ylabel={$\| u - u_h \|_{\infty} / \| u \|_{\infty}$},
        ymode=log,
        legend entries={$m=2$,$m=4$,$m=8$,$m=16$,$m=32$,$m=64$},
        legend pos=south west,
        x=1.2cm,y=0.8cm,
        ymin=0.0001, ymax=0.3,
    ]
    \addplot plot coordinates { (0, 0.157935) (1, 0.093944) (2, 0.057234) (3, 0.029642) (4, 0.026702) (5, 0.025229) (6, 0.024303) };
    \addplot plot coordinates { (0, 0.157935) (1, 0.083140) (2, 0.032728) (3, 0.012080) (4, 0.009543) (5, 0.007900) (6, 0.006878) };
    \addplot plot coordinates { (0, 0.178669) (1, 0.087210) (2, 0.034331) (3, 0.005554) (4, 0.003409) (5, 0.002464) (6, 0.002142) };
    \addplot plot coordinates { (0, 0.254317) (1, 0.186331) (2, 0.015656) (3, 0.004477) (4, 0.001677) (5, 0.000962) (6, 0.000728) };
    \addplot plot coordinates { (0, 0.255076) (1, 0.229091) (2, 0.042021) (3, 0.004912) (4, 0.001540) (5, 0.000503) (6, 0.000258) };
    \addplot plot coordinates { (0, 0.255076) (1, 0.229091) (2, 0.016613) (3, 0.004954) (4, 0.001470) (5, 0.000503) (6, 0.000135) };
    \end{axis}
\end{tikzpicture}
\end{center}
\caption{Relative $L^\infty(\omega)$-error for the test problem of Experiment \ref{exp:Bakelman}.}
\label{fig:Bakelman:rates}
\end{figure}

\begin{experiment}[Non-convex boundary data on convex domain] \label{exp:Bakelman}
Now
\[
\Omega = B_1(0) \cup \bigr[ (0,1) \times (0,1) \bigr] = \{ x \in \bR^2 : \| x \| < 1 \} \cup \bigr[ (0,1) \times (0,1) \bigr]
\]
is a convex but not strictly convex domain. On the vertical straight boundary segment $\Gamma := \partial \Omega \cap \bigl[ \{ 1 \} \times (0,1) \bigr]$ we impose the non-convex boundary data $g(x) = x_2 (1-x_2)$ while on the remainder $\partial \Omega \setminus \Gamma$ we set $g(x) = x_1^4 - 1$. We select $f = 0$ on $\Omega$, implying that the graph of the solution $u$ is a surface of vanishing Gauss curvature.

A numerical solution is shown in Figure \ref{fig:Bakelman:plot}. In the vicinity of $\Gamma$ on the far side of the plot the numerical solution is nearly vertical, interpolating on $\Gamma$ the data $g(x) = x_1 (1 - x_1)$ while attaining at the interior nodes neighbouring $\Gamma$ approximately the value \mbox{$x_1^4 - 1 \approx 0$}.

Recalling from \eqref{eq:equivalence} that the Bellman formulation enforces convexity on the domain $\Omega$ of the differential operator, it appears that the boundary data is only extended into the interior of the domain in as far as convexity permits. Indeed, according Figure \ref{fig:Bakelman:rates} the numerical solutions converge on the subdomain $\omega := \{ x \in \Omega : x_1 < 0.95 \}$ under mesh refinement to $u(x) = x_1^4 - 1$; we see a relative $L^\infty(\omega)$ error of $1.35 \cdot 10^{-4}$ on a mesh with $304129$ DoFs and $m = 64$.

The table gives the numerical values of the error for the best choices of $m$ on a given mesh, and the associated number of Newton iterations:

\begin{center}
\begin{tabular}{r||r|r|r|r|r|r|r}
DoFs & 145 & 329 & 1249 & 4865 & 19201 & 76289 & 304129\\ \hline
$m$ & 2 & 4 & 16 & 16 & 64 & 64 & 64\\ 
rel.~$L^\infty$-error & 1.6e-1 & 8.3e-2 & 1.6e-2 & 4.5e-3 & 1.5-3 & 5.0e-04 & 1.4e-4\\
\# Newton & 6 & 6 & 9 & 10 & 12 & 14 & 19
\end{tabular}
\end{center}

Overall the largest number of Newton iterations in this computational experiment is $19$, which occurs on the finest mesh with $m = 64$.
\end{experiment}

To give an interpretation of Experiment \ref{exp:Bakelman} we review a result due to Bakelman \cite{B}. Let $v$ be a bounded convex function on $\Omega$ and $\overline{Co}(\Gamma_v)$ be the closed convex hull of the graph $\Gamma_v \subset \bR^{2+1}$ of $v$. 
Then the function
\[
\nu_v : \; \partial \Omega \to \bR, \; x \mapsto \inf \{ h \in \bR : \, (x,h) \in \overline{Co}(\Gamma_v) \}
\]
is called the border of the function $u$. We denote by $\mathcal L$ the set of functions $v$ which satisfy the Monge-Amp\`ere differential equation in the Aleksandrov sense and
\[
\nu_v(x) \le g(x) \qquad  \forall x \in \partial \Omega.
\]
If the domain $\Omega$ is strictly convex then Aleksandrov and viscosity solutions coincide for the simple Monge-Amp\`ere equation: See \cite{G01} for the proof when the continuous $f$ is positive. This argument can be extended to the case of non-negative $f$ with the tools of \cite{dPF}---alternatively one can use the uniqueness of solutions and \cite{FJ}.

In the context of the Gauss curvature problem, Bakelman shows under a data condition that on bounded strictly convex domains $\Omega$ the set $\mathcal L$ is non-empty. Moreover, there is a unique $u \in \mathcal L$ such that $\nu_u \ge \nu_v$ for all $v \in \mathcal L$. This $u$ is considered to be the generalised solution of the boundary value problem.

The numerical solutions of Experiment \ref{exp:Bakelman} evidently converge to $u(x) = x_1^4 - 1$ in $L^\infty_{\rm loc}(\Omega)$, so on relatively compact subsets of $\Omega$. The border of $u$ is its trace. As the numerical solution attains the continuum of values between $g(x)$ and $\nu_u(x) = u(x)$ in the vicinity of any $x \in \partial \Omega$, the visual impression is that of a convergence to a multi-valued limit which attains the interval $[\nu_u(x), g(x)]$ at $x$. It is therefore appealing to adopt this multi-valued interpretation for the purposes of this text.

Starting from \cite{B}, Froese recently studied a wide-stencil method for the Gauss curvature problem with non-classical boundary conditions on convex domains \cite{F}. There, besides the definition with the border function $\nu_u$, also viscosity boundary conditions as in \cite{UG} are discussed.

\section{The simple Monge--Amp\`ere equation on non-convex domains without a classical subsolution}

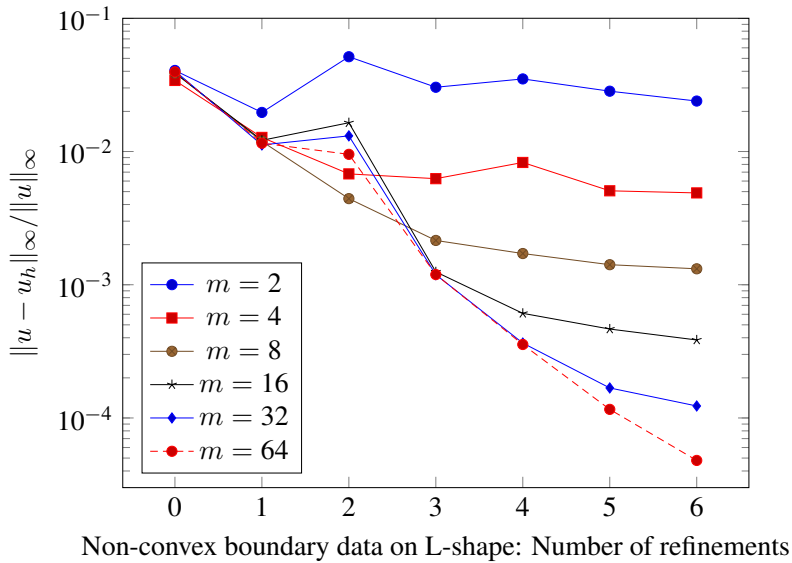
\begin{figure}[t]
\begin{center}
\begin{tikzpicture}
    \begin{axis}[
        xlabel={Non-convex boundary data on L-shape: Number of refinements},
        ylabel={$\| u - u_h \|_{\infty} / \| u \|_{\infty}$},
        ymode=log,
        legend entries={$m=2$,$m=4$,$m=8$,$m=16$,$m=32$,$m=64$},
        legend pos=south west,
        x=1.2cm,y=0.8cm,
        ymin=0.00003, ymax=0.1,
    ]

    \addplot plot coordinates { (0, 0.040580) (1, 0.019647) (2, 0.051444) (3, 0.030356) (4, 0.035057) (5, 0.028360) (6, 0.023919) };
    \addplot plot coordinates { (0, 0.034143) (1, 0.012758) (2, 0.006792) (3, 0.006253) (4, 0.008275) (5, 0.005070) (6, 0.004889) };
    \addplot plot coordinates { (0, 0.038878) (1, 0.011910) (2, 0.004429) (3, 0.002150) (4, 0.001715) (5, 0.001414) (6, 0.001317) };
    \addplot plot coordinates { (0, 0.038877) (1, 0.012129) (2, 0.016437) (3, 0.001249) (4, 0.000610) (5, 0.000465) (6, 0.000385) };
    \addplot plot coordinates { (0, 0.040269) (1, 0.011203) (2, 0.013103) (3, 0.001190) (4, 0.000365) (5, 0.000168) (6, 0.000123) };
    \addplot plot coordinates { (0, 0.039595) (1, 0.011546) (2, 0.009522) (3, 0.001192) (4, 0.000356) (5, 0.000116) (6, 0.000048) };
    \end{axis}
\end{tikzpicture}
\end{center}
\caption{Relative $L^\infty(\omega)$-error for the test problem of Experiment \ref{exp:nonconvexLshape}.}
\label{fig:nonconvexLshape:rates}
\end{figure}

\begin{figure}[p]
\begin{center}
\includegraphics[width=12.5cm]{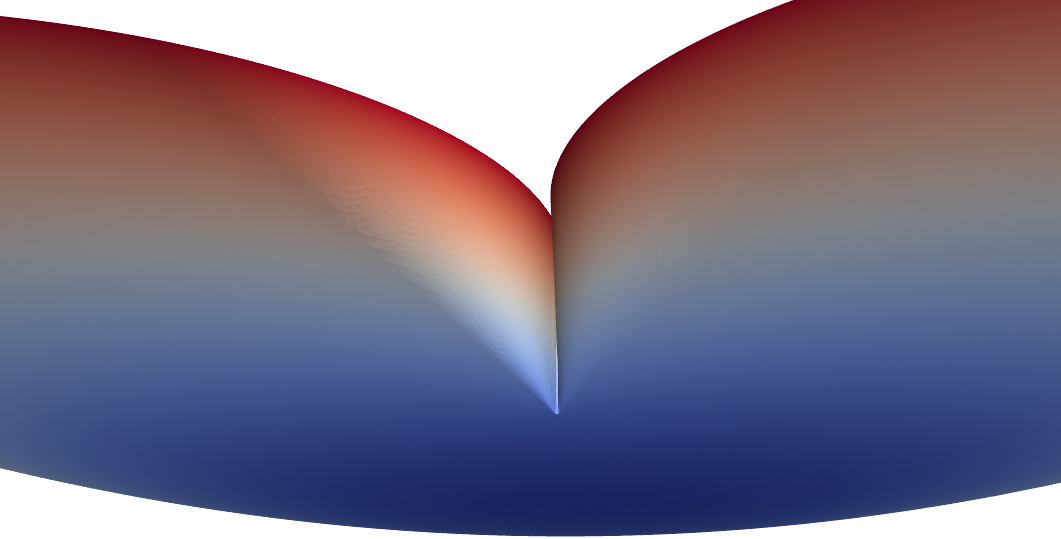}\\
(a) Close-up view onto the numerical solution at non-convex part of $\Omega$.

\bigskip

\includegraphics[width=12cm]{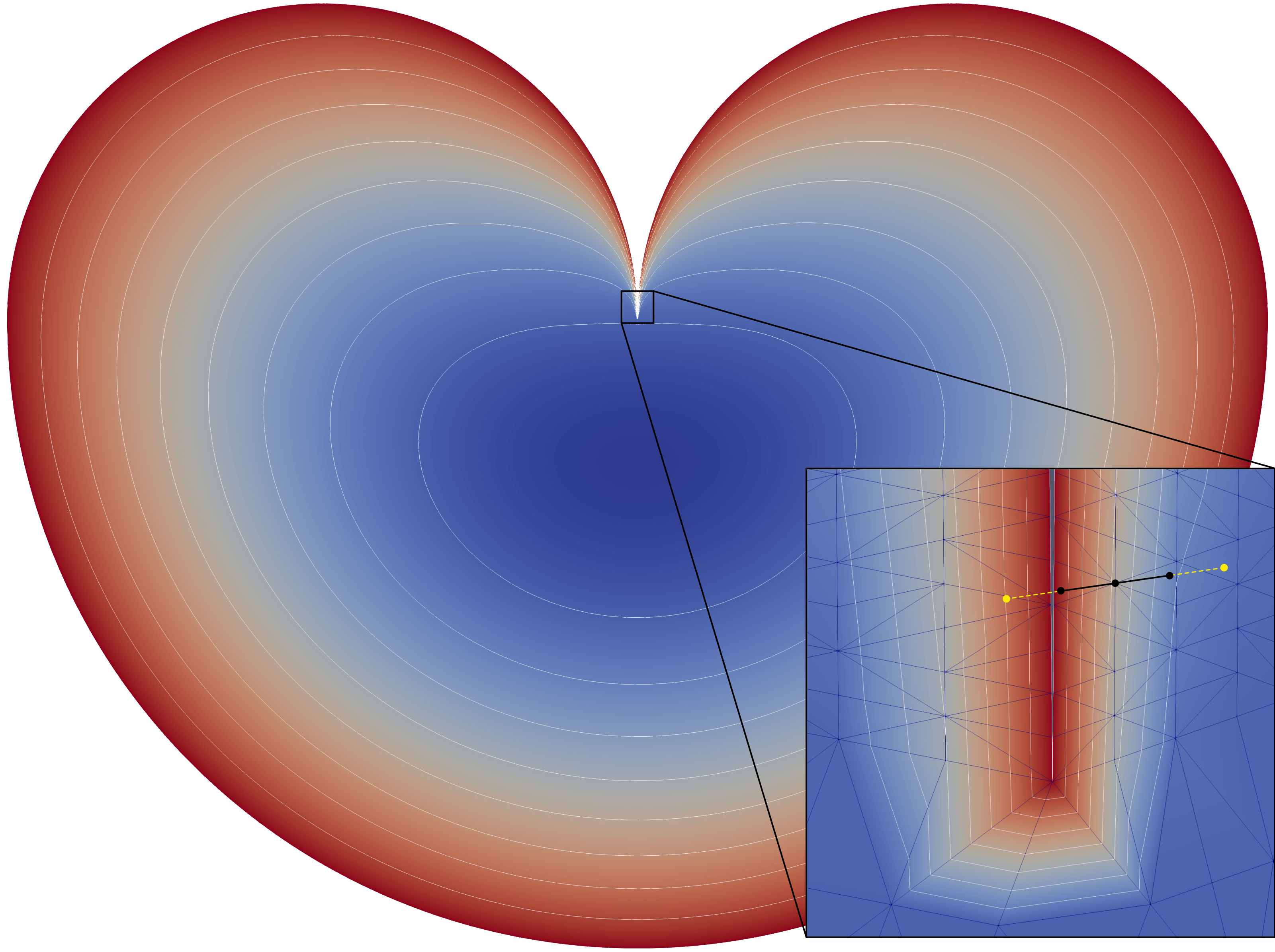}\\
(b) Numerical solution with contour lines.
\end{center}
\caption{It appears that the numerical solution of Experiment \ref{exp:heart} approximates the transition from a multi-valued boundary condition at the origin to a classical boundary condition on the remainder of the boundary; see the close-up (a) and also in the (b) the concentration of contour lines within the width of a single element at the origin. A cartoon stencil in (b) shows that careful scaling of the stencil is important to ensure that the finite differences approximate the local features of the solution.}
\label{fig:heart:plot}
\end{figure}

%

As far as we are aware there is no systematic analysis in the literature which extends \cite{B,GS93,G98} to boundary value problems where non-convex boundary data is imposed on a non-convex domain. 

Yet, also in this setting the numerical scheme of \cite{FJ} is guaranteed to have a unique solution for any mesh. From our point of view this makes numerical experiments interesting. 

Indeed, equations of Monge--Amp\`ere type have been proposed for physical and biological models where the application does not justify to impose convexity on $\Omega$ or $g$. An example is the Monge--Amp\`ere Keller--Segel system of chemotaxis \cite{B08,HL}. There the solution of a Monge--Amp\`ere equation describes the density of a chemical substance, whose physical domain must not necessarily be convex.

\begin{experiment}[Non-convex boundary data on the L-shape] \label{exp:nonconvexLshape}
Combining in spirit the domain of Experiment \ref{exp:quarticLshape} with the boundary data of Experiment \ref{exp:Bakelman}, we set $g(x) = x_2 (1-x_2)$ on left boundary segment $\Gamma := \partial \Omega \cap \bigl[ \{ -1 \} \times (0,1) \bigr]$ of the L-shape $\Omega = \bigl[ (0,1) \times (-1,1) \bigr] \cup \bigl[ (-1,1) \times (0,1) \bigr]$. On the remainder $\partial \Omega \setminus \Gamma$ we set $g(x) = x_1^4 - 1$ and on $\Omega$ we select $f = 0$.

Similarly to the previous experiment, numerical solutions approximate a multi-valued solution in the vicinity of $\Gamma$, while on the subdomain $\omega := \{ x \in \Omega : x_1 > -0.95 \}$ the numerical solutions converge under mesh refinement to $u(x) = x_1^4 - 1$. We see a relative $L^\infty(\omega)$ error of $4.8 \cdot 10^{-5}$ on a mesh with $168961$ DoFs and $m = 64$, cf. Figure \ref{fig:nonconvexLshape:rates}. Overall the largest number of Newton iterations in this computational experiment is $25$, which occurs on the finest mesh with $m = 32$.
\end{experiment}

The multi-valued behaviour may also be exhibited by problems with convex or even vanishing boundary data, when no smooth classical subsolution exists.

\begin{figure}[t]
\begin{center}
\includegraphics[width=6cm]{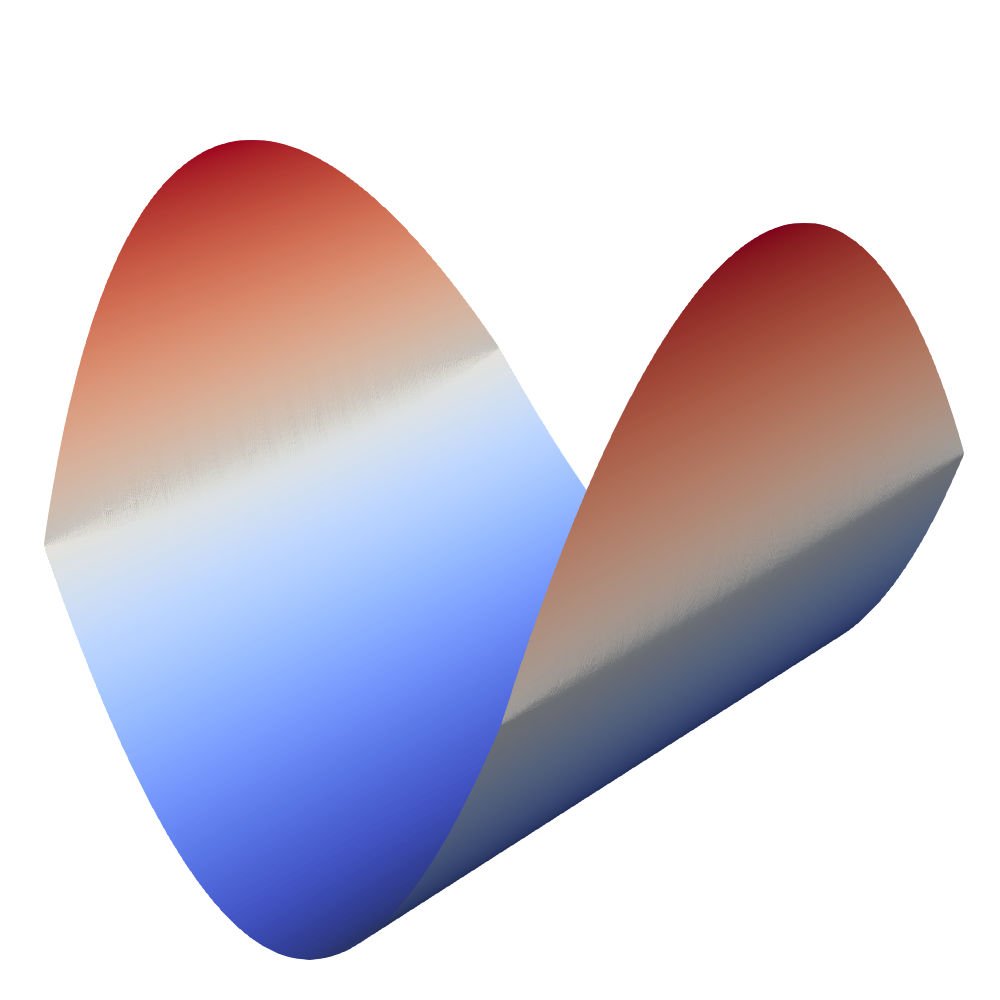} \qquad \includegraphics[width=6cm]{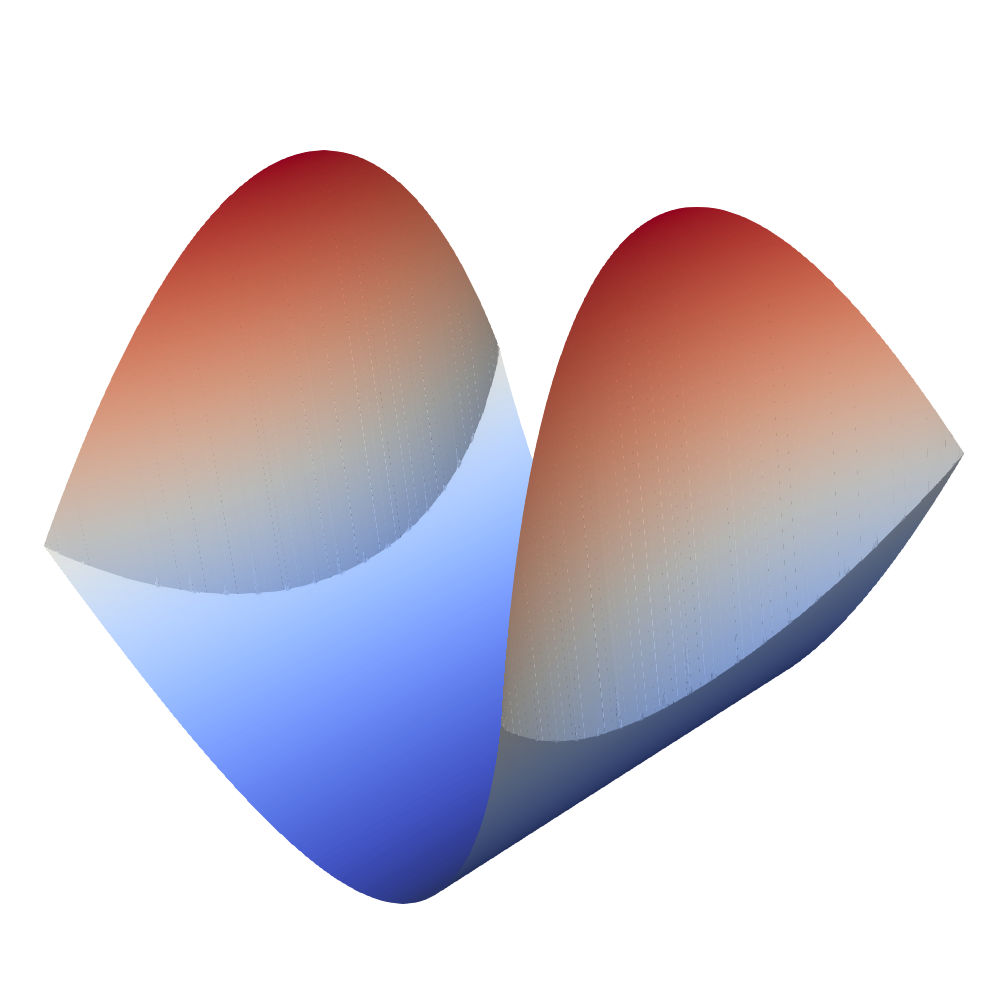}\\
(a) \hspace{7cm} (b)
\end{center}
\caption{The figure shows the numerical solutions in the two domains with $c = 3$ of Experiment \ref{exp:bent-square}, on the left with convex and on the right with concave faces.}
\label{fig:bent-square:plot}
\end{figure}

\begin{experiment}[H\"older domain] \label{exp:heart}
We now consider the heart-shaped domain
\begin{align*}
\Omega = & \, \bigl[ \bigl( B_{0.5}((0,0.5)) \cup B_{0.5}((0,-0.5)) \bigr) \cap \{ x \in \bR^2 : x_1 \ge 0 \} \bigr]\\
& \, \cup \bigl[ B_1((0,0)) \cap \{ x \in \bR^2 : x_1 \le 0 \} \bigr],
\end{align*}
where again $B_\delta(x)$ denotes the ball with radius $\delta$ and centre $x$. This domain does not satisfy a Lipschitz condition at the origin, as can also be seen in Figure \ref{fig:heart:plot} (b). In this example the boundary data $g = 0$ vanishes, while $f = 1$.

The plots of Figure \ref{fig:heart:plot} give the impression that at the origin the numerical solution approximates a multi-valued function but then transitions to a classical single-valued boundary condition. The Newton method stopped the computation with 166,176 DoFs after 16 iterations, when a Newton step size of $5 \cdot 10^{-8}$ was reached.

A note about the stencil diameter: As the cartoon stencil in Figure \ref{fig:heart:plot} (b) illustrates, a too simple implementation a wide second-order central difference with the centre on the right of the re-entrant boundary could have a node in $\Omega$ which is on the left of the re-entrant boundary. On non-convex domains stencils should be scaled so that not just the nodes of the stencil belong to $\Omega$, but also any points between them.
\end{experiment}

\begin{experiment}[Bent square] \label{exp:bent-square}
In this last experiment we solve the degenerate equation on a sequence of curved polyhedra, starting from a uniformly convex geometry and leading to polyhedra with concave faces:

\begin{center}
\includegraphics[width=1.8cm]{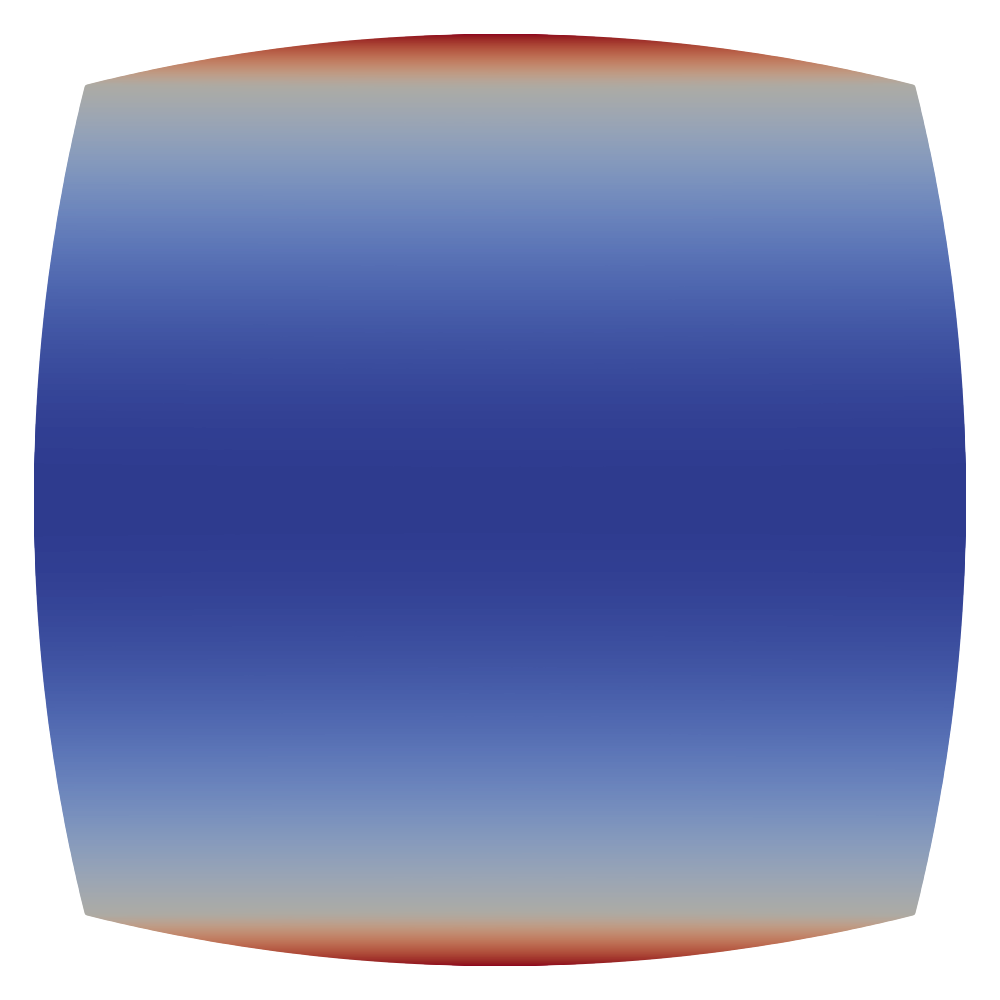}  \includegraphics[width=1.8cm]{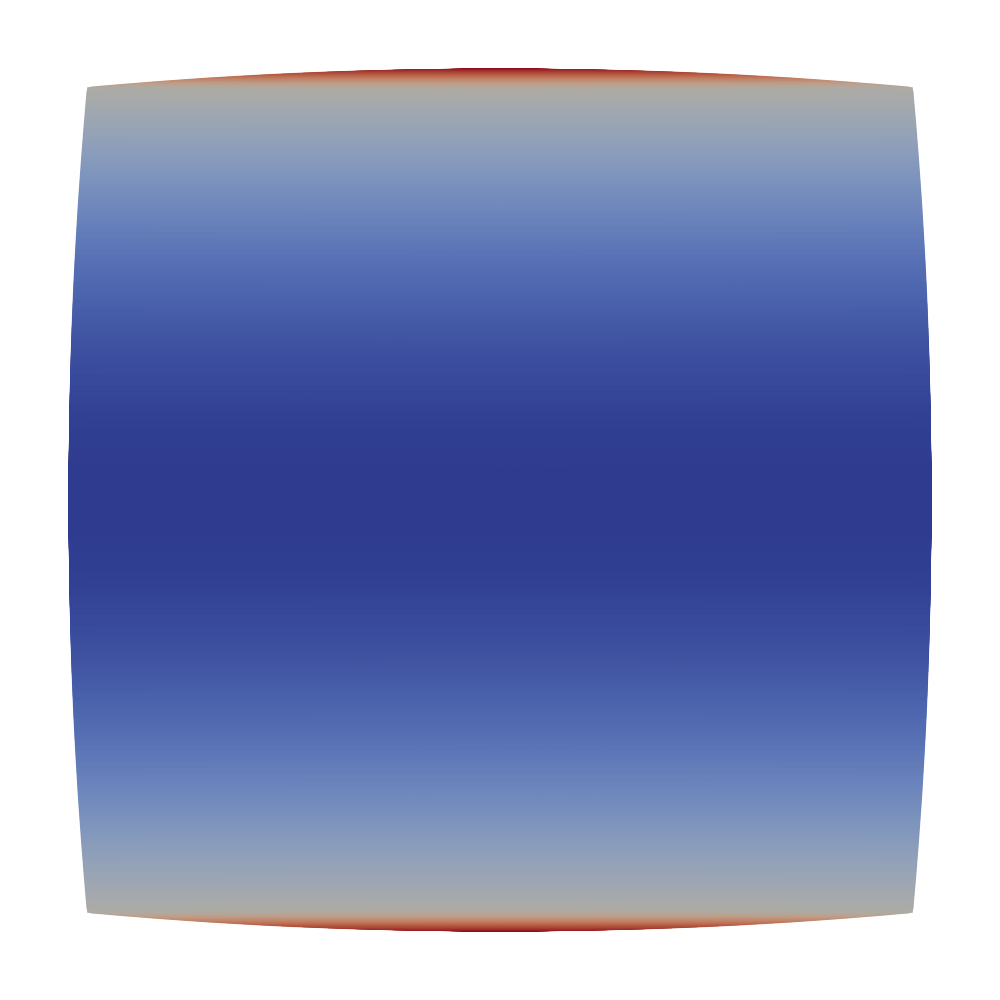}  \includegraphics[width=1.8cm]{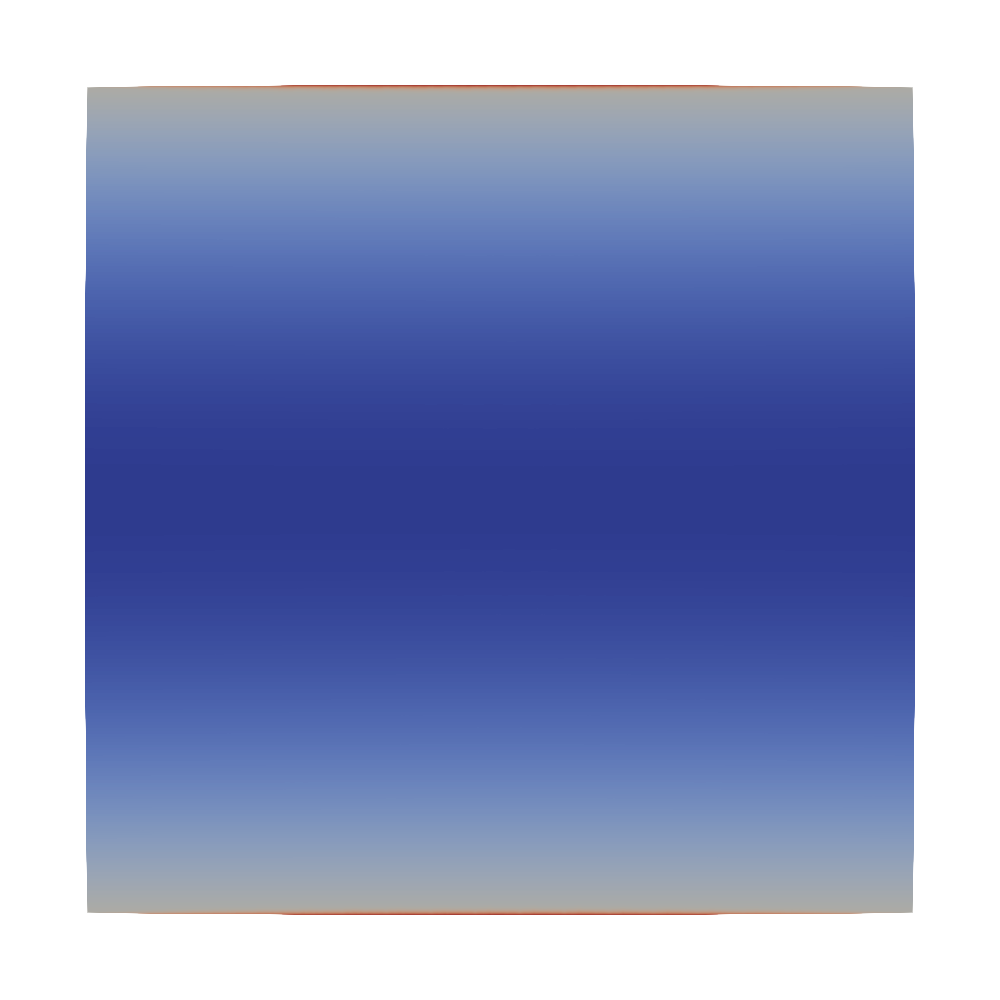}  \includegraphics[width=1.8cm]{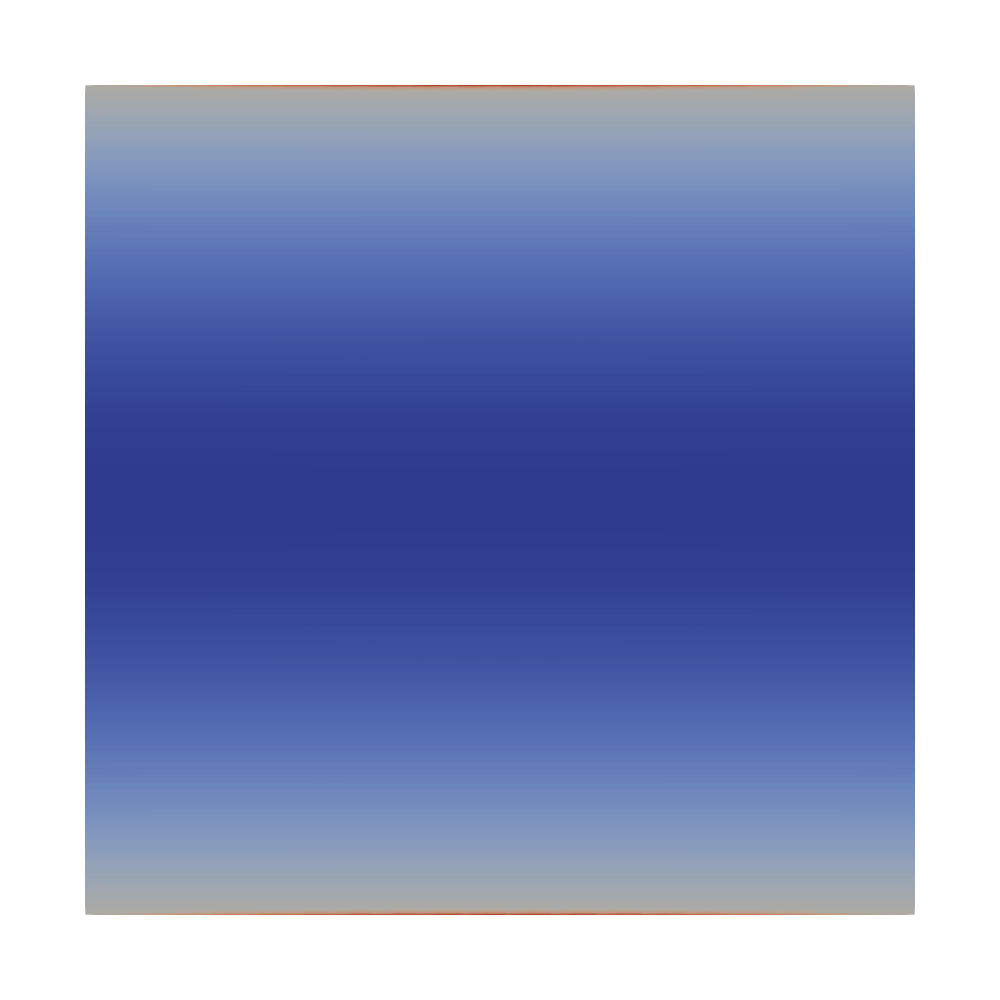}  \includegraphics[width=1.8cm]{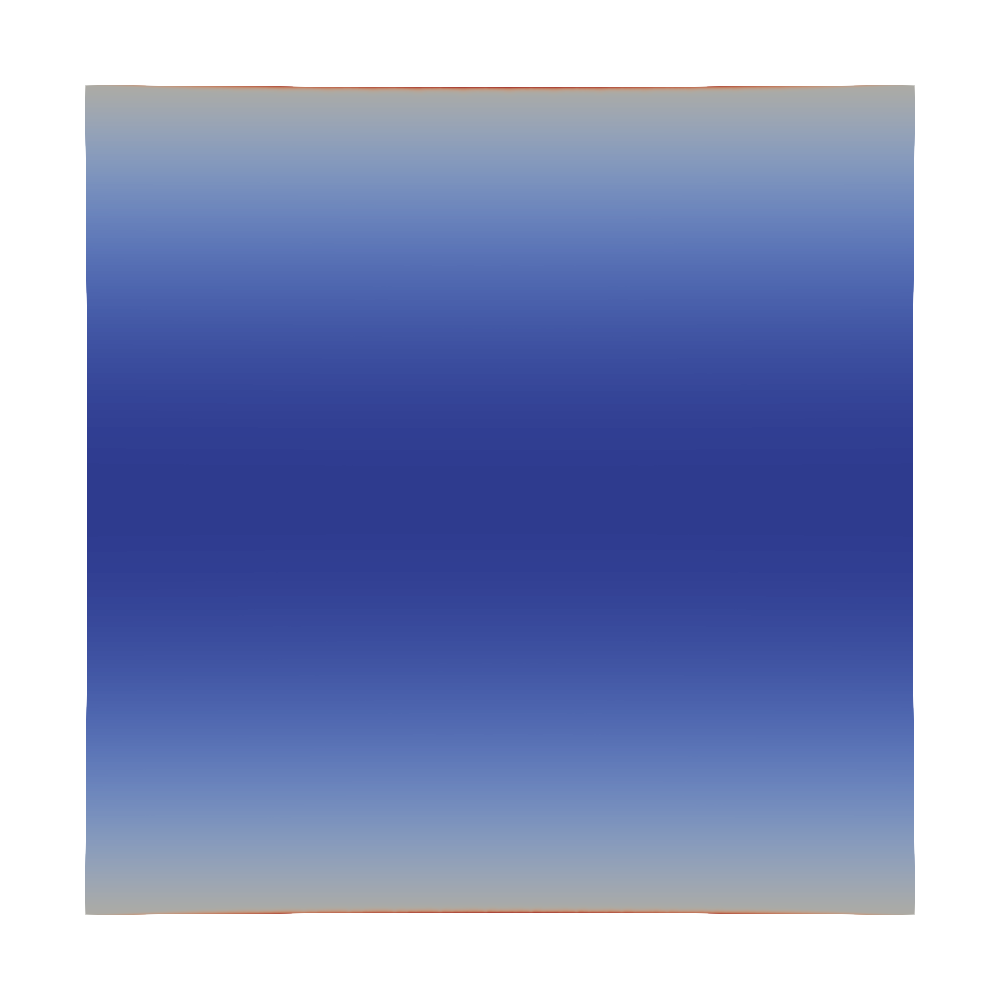}  \includegraphics[width=1.8cm]{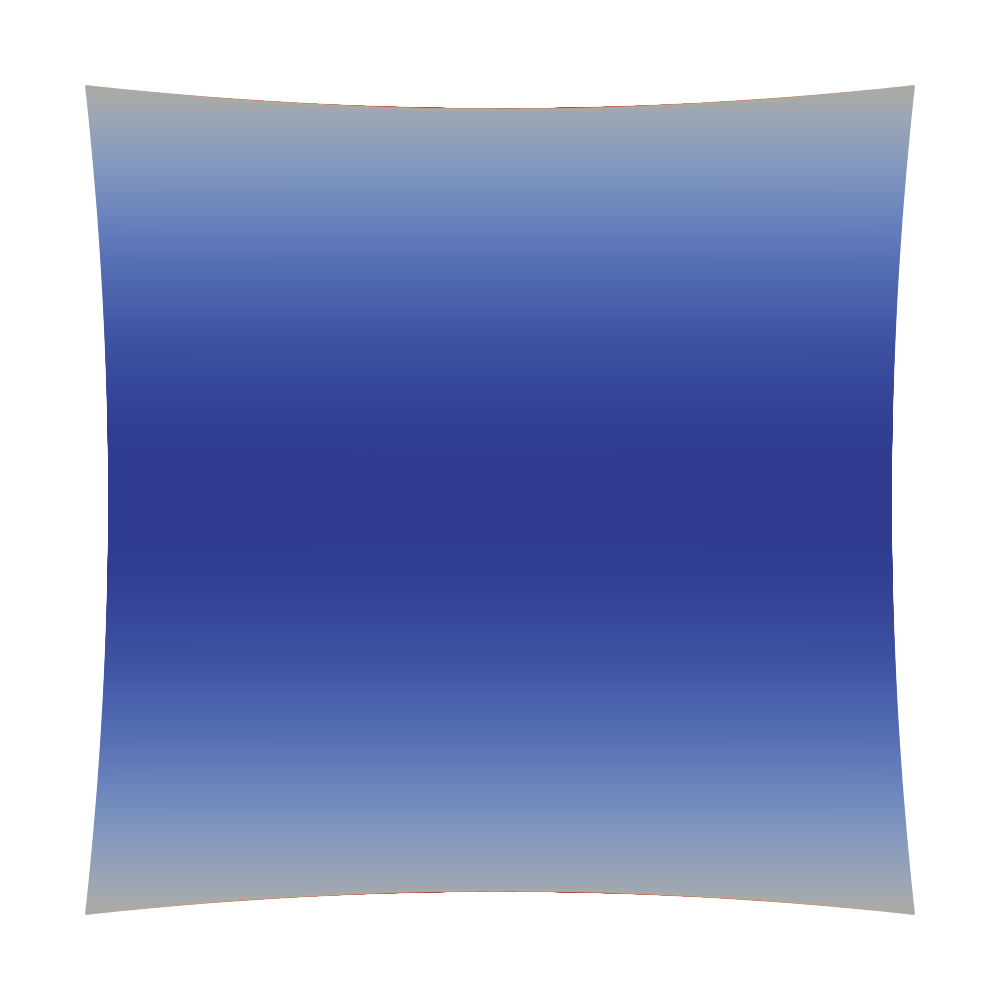}  \includegraphics[width=1.8cm]{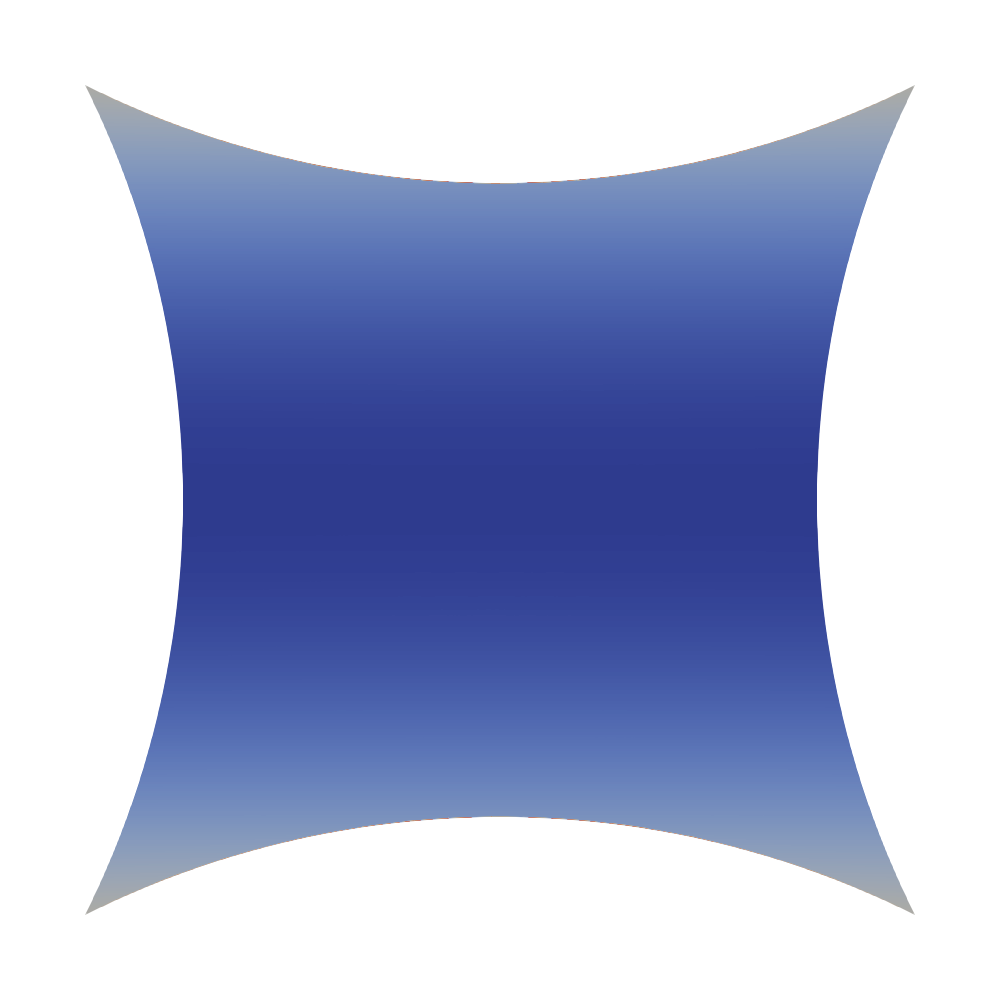}
\end{center}

To be precise: The first three domains are formed from the intersection of four circles with the midpoints $(c, 0)$, $(0, c)$, $(-c, 0)$, $(0, -c)$ where $c$ is $3$, $10$, $100$ respectively and where the radii are chosen so that the circles intersect at the points $(1,1)$, $(-1,1)$, $(1,-1)$ and $(-1,-1)$. The fourth domain is the square centred at the origin with side length $2$. The remaining three domains are obtained by subtracting from this square the points which lie in circles with the same centres and intersection points as above, but smaller radii.

With $f = 0$ the graphs of the exact solutions are (when restricted to the interior) surfaces of vanishing Gauss curvature. The function $g$ is equal to the convex function $x_2^2 - 1$ where $x_1 > x_2$ and equal to the concave function $1 - x_1^2$ where $x_1 \le x_2$. We note that $g$ is continuous on $\partial \Omega$. Two of the resulting solutions are plotted in Figure \ref{fig:bent-square:plot}. In the interior the numerical solutions do not significantly vary from the function $x_2^2 - 1$, which already defined the convex part of the boundary data. This is notable because on strictly convex domains the boundary function $g$ is attained on all of $\partial \Omega$ in the classical single-valued sense, while on the other domains the multi-valued cut-off mechanisms appears to screen out concave sections of $g$. Indeed on $\omega = (-3/4, 3/4)^2$ we find the following relative $L^\infty$ errors, taking $x_2^2 - 1$ as the reference solution:

\begin{center}
\begin{tabular}{r||r|r|r|r|r|r|r}
domain & 1 & 2 & 3 & 4 & 5 & 6 & 7\\ \hline
rel.~$L^\infty(\omega)$-error & 1.4e-3 & 1.3e-3 & 1.4e-3 & 1.4e-3 & 1.2e-3 & 1.2e-3 & 8.7e-4\\
\# Newton & 10 & 10 & 9 & 9 & 14 & 20 & 20
\end{tabular}
\end{center}

The table also shows an increase in the required Newton iterations to achieve a Newton step size of $5 \cdot 10^{-8}$ when the faces of $\Omega$ are concave. The number of DoFs vary in this experiment between 108,124 and 140,324, depending on the domain.
\end{experiment}

\end{document}